# CRITICAL RANDOM GRAPHS: DIAMETER AND MIXING TIME[1]

By Asaf Nachmias and Yuval Peres

*University of California, Berkeley and Microsoft Research*

Let $\mathcal{C}_1$ denote the largest connected component of the critical Erdős–Rényi random graph $G(n, \frac{1}{n})$. We show that, typically, the diameter of $\mathcal{C}_1$ is of order $n^{1/3}$ and the mixing time of the lazy simple random walk on $\mathcal{C}_1$ is of order $n$. The latter answers a question of Benjamini, Kozma and Wormald. These results extend to clusters of size $n^{2/3}$ of $p$-bond percolation on any $d$-regular $n$-vertex graph where such clusters exist, provided that $p(d-1) \leq 1 + O(n^{-1/3})$.

**1. Introduction.** The Erdős–Rényi random graph $G(n, p)$ is obtained from the complete graph on $n$ vertices by retaining each edge with probability $p$ and deleting it with probability $1 - p$, independently of all other edges. Fountoulakis and Reed [12] and Benjamini, Kozma and Wormald [5] proved that the mixing time of a random walk on the largest connected component $\mathcal{C}_1$ of $G(n, \frac{\theta}{n})$ with $\theta > 1$ is of order $\log^2(n)$ with high probability. The latter authors asked what the mixing time is in the critical random graph $G(n, 1/n)$. The next theorem (a special case of our main result, Theorem 1.2) answers their question and also determines the diameter of $\mathcal{C}_1$ in this case.

*Terminology.* A *lazy simple random walk* on a graph $G = (V, \mathcal{E})$ is a Markov chain on $V$ with transition probabilities $\mathbf{p}(x, y) = \frac{1}{2 \deg(x)}$ if $(x, y) \in \mathcal{E}$ and $\mathbf{p}(x, x) = \frac{1}{2}$ for all $x \in V$. It has stationary distribution $\pi$ given by $\pi(x) = \frac{\deg(x)}{2|\mathcal{E}|}$. The *mixing time* of the lazy random walk on $G$ is

$$T_{\mathrm{mix}}(G) = T_{\mathrm{mix}}(G, 1/4) = \min\{t : \|\mathbf{p}^t(x, \cdot) - \pi(\cdot)\|_{\mathrm{TV}} \leq 1/4, \text{ for all } x \in V\},$$

where $\|\mu - \nu\|_{\mathrm{TV}} = \max_{A \subset V} |\mu(A) - \nu(A)|$ is the total variation distance.

Received February 2007; revised July 2007.
[1]Supported in part by NSF Grants DMS-02-44479 and DMS-01-04073.

THEOREM 1.1. *Let $\mathcal{C}_1$ denote the largest connected component of $G(n, \frac{1+\lambda n^{-1/3}}{n})$ for $\lambda \in \mathbb{R}$. Then, for any $\varepsilon > 0$, there exists $A = A(\varepsilon, \lambda) < \infty$ such that for all large $n$:*

- $\mathbf{P}(\mathrm{diam}(\mathcal{C}_1) \notin [A^{-1} n^{1/3}, A n^{1/3}]) < \varepsilon$;
- $\mathbf{P}(T_{\mathrm{mix}}(\mathcal{C}_1) \notin [A^{-1} n, A n]) < \varepsilon$.

For $\theta > 1$, the diameter of the largest component in $G(n, \frac{\theta}{n})$ is typically of order $\log n$; see [6, 9, 11]. For $\theta < 1$, the diameter of the largest component in $G(n, \frac{\theta}{n})$ is typically of order $\sqrt{\log n}$, but there are components of smaller cardinality with diameter of order $\log n$; see [17]. In $G(n, \frac{1}{n})$, it is natural to expect that the diameter of $\mathcal{C}_1$ will be of order $n^{1/3}$ since a random tree on $m$ vertices typically has diameter of order $\sqrt{m}$ (see, e.g., [14]) and with probability bounded below, $\mathcal{C}_1$ is a tree with roughly $n^{2/3}$ vertices. Indeed, if $\mathcal{C}_1$ is a tree, then it is easy to deduce the bounds on the diameter and the upper bound on the mixing time. However, the probability that $\mathcal{C}_1$ is a tree does not tend to 1 as $n \to \infty$.

We state our main result in the more general setting of *percolation on finite graphs*. Given a finite graph $G$ and $p \in (0, 1)$, the random subgraph $G_p$ is obtained from bond percolation with parameter $p$, that is, each edge of $G$ is (independently) retained with probability $p$ and erased with probability $1 - p$. The next theorem states that when $G$ has maximum degree at most $d \in [3, n-1]$ and $p \leq \frac{1+\lambda n^{-1/3}}{d-1}$, if $G_p$ typically has components of order $n^{2/3}$, then, with high probability, all such components will have diameter of order $n^{1/3}$ and mixing time of order $n$. The components are unlikely to be larger than about $n^{2/3}$, by part (b) of the theorem.

THEOREM 1.2. *Let $G$ be a graph on $n$ vertices with maximum degree at most $d \in [3, n-1]$. For $0 < p < 1$, denote by $\mathbf{CO}(G_p)$ the collection of connected components of the percolation subgraph $G_p$. For $\mathcal{C} \in \mathbf{CO}(G_p)$, let $\mathcal{E}(\mathcal{C})$ denote the edge set in $\mathcal{C}$ and recall that $T_{\mathrm{mix}}(\mathcal{C})$ is the mixing time of lazy simple random walk on $\mathcal{C}$. If $p \leq \frac{1+\lambda n^{-1/3}}{d-1}$ for some fixed $\lambda \in \mathbb{R}$, then for any $\varepsilon > 0$ and $\beta > 0$, there exists $A = A(\varepsilon, \beta, \lambda) < \infty$ such that for all large $n$:*

(a) $\mathbf{P}(\exists \mathcal{C} \in \mathbf{CO}(G_p) \text{ with } |\mathcal{C}| > \beta n^{2/3} \text{ and } \mathrm{diam}(\mathcal{C}) \notin [A^{-1} n^{1/3}, A n^{1/3}]) < \varepsilon$;

(b) $\mathbf{P}(\exists \mathcal{C} \in \mathbf{CO}(G_p) \text{ with } |\mathcal{E}(\mathcal{C})| > A n^{2/3}) < \varepsilon$;

(c) $\mathbf{P}(\exists \mathcal{C} \in \mathbf{CO}(G_p) \text{ with } |\mathcal{C}| > \beta n^{2/3} \text{ and } T_{\mathrm{mix}}(\mathcal{C}) \notin [A^{-1} n, A n]) < \varepsilon$.

Since $G(n, p)$ is $G_p$, where $G$ is the complete graph on $n$ vertices, Theorem 1.1 is an immediate consequence of Theorem 1.2 and the following fact, first proven in [10] (see also [1, 18, 21]): in $G(n, \frac{1+\lambda n^{-1/3}}{n})$,

$$\liminf_n \mathbf{P}(|\mathcal{C}_1| > \beta n^{2/3}) \to 1 \qquad \text{as } \beta \to 0.$$

Similarly, one can deduce the same results for the $j$th largest component of the random graph $G(n, \frac{1+\lambda n^{-1/3}}{n})$ for any constant $j$. Other examples of $d$-regular graphs $G$ where the hypothesis of Theorem 1.2 are satisfied (and for $p = \frac{1+\lambda n^{-1/3}}{d-1}$, there are components of size greater than $\beta n^{2/3}$ in $G_p$ with probability bounded away from 0) are uniform random $d$-regular graphs (see [22]) and the Cartesian square of a complete graph (see [13] and Theorem 1.3 of [7]).

We also show that the maximal diameter over all components is typically at most $O(n^{1/3})$ and only components with cardinality of order $n^{2/3}$ can achieve this diameter. This contrasts with the subcritical case where, as noted above, there are components with diameter of order $\log n$, but the diameter of the largest component is typically of order $\sqrt{\log n}$; see [17].

THEOREM 1.3. *Under the conditions of Theorem 1.2, we have that for any $\varepsilon > 0$ and $\beta > 0$, there exists $A = A(\varepsilon, \beta, \lambda) < \infty$ such that*

(1.1) $$\mathbf{P}(\exists \mathcal{C} \in \mathbf{CO}(G_p) \text{ with } \operatorname{diam}(\mathcal{C}) > An^{1/3}) < \varepsilon.$$

*Furthermore, for any $D_1 > 0$, there exists $D_2 > 0$ such that for all $M < n^{2/3}/2$,*

(1.2) $$\mathbf{P}(\exists \mathcal{C} \in \mathbf{CO}(G_p): |\mathcal{C}| < M \text{ and } \operatorname{diam}(\mathcal{C}) > D_2\sqrt{M \log(n/M^{3/2})}) \leq \left(\frac{M^{3/2}}{n}\right)^{D_1}.$$

For the random graph $G(n,p)$ and the random $d$-regular graph, we can prove a stronger tail bound on the diameter of the connected components.

PROPOSITION 1.4. *Assume that $G$ is the complete graph $K_n$. If $p \leq \frac{1+\lambda n^{-1/3}}{n}$, then there exists $c = c(\lambda) > 0$ such that*

$$\mathbf{P}(\exists \mathcal{C} \in \mathbf{CO}(G_p) \text{ with } \operatorname{diam}(\mathcal{C}) > An^{1/3}) \leq e^{-cA^{3/2}}.$$

*The same inequality [with $c = c(d, \lambda)$] holds if $G$ is a random uniform $d$-regular graph on $n$ vertices where $d$ is fixed and $p \leq \frac{1+\lambda n^{-1/3}}{d-1}$.*

The rest of the paper is organized as follows. In the remainder of this section, we discuss the intuition and idea for the proofs of this paper. In Section 2, we present some preliminaries. For ease of exposition, we first prove Theorem 1.2 for the case $\lambda \leq 0$ (i.e., when $p \leq \frac{1}{d-1}$) and defer the case of $\lambda > 0$ to Section 6. Parts (a) and (b) of Theorem 1.2 are established in Section 3 under the assumption $\lambda \leq 0$, which allows for a very short proof. The upper bound on the mixing time follows easily, so we present it in Section 4. The lower bound on the mixing time is given in Section 5. In Section 6, we prove Theorem 1.2 in its full generality and also Theorem 1.3 and Proposition 1.4. These proofs all have a common element, Lemma 6.2, which shows that for the parameters we are considering, the diameter of a component $\mathcal{C}$ is unlikely to be much larger than the square root of the number of vertices in $\mathcal{C}$.

*Proofs idea.* We start by understanding the reason why, in $p$-bond percolation on any $d$-regular graph with $p \leq \frac{1}{d-1}$, the largest component will be of size no more than $n^{2/3}$. This is stated in part (b) of Theorem 1.2; here, we present a slightly different approach for the purpose of exposition.

Fix a vertex $v$ of the graph and let $\mathcal{C}(v)$ be the component of $G_p$ containing $v$. Let $T$ be the surviving population of a Galton–Watson branching process in which the initial particle (the root) has progeny distribution $\mathrm{Bin}(d,p)$ and all other particles have progeny distributions $\mathrm{Bin}(d-1,p)$. Observe that since $G$ is $d$-regular, we can naturally couple such that $|\mathcal{C}(v)| \leq |T|$. It is a well-known fact in the theory of branching processes (see [14]) that as long as $p \leq \frac{1}{d-1}$, there is a constant $c > 0$ such that for all $M > 0$, we have

$$\mathbf{P}(|T| \geq M) \leq \frac{c}{\sqrt{M}}. \tag{1.3}$$

In fact, $|T|$ is the first time a certain random walk with mean 0 and bounded increments visits 0 and (1.3) is also classical in the theory of random walks. By our coupling, we have the same tail upper bound (1.3) for $|\mathcal{C}(v)|$. Let $X = |\{v : |\mathcal{C}(v)| \geq M\}|$ so that $\mathbb{E}X \leq \frac{cn}{\sqrt{M}}$. Observe that if the largest component $\mathcal{C}_1$ has at least $M$ vertices, then $X \geq M$. Therefore,

$$\mathbf{P}(|\mathcal{C}_1| \geq M) \leq \mathbf{P}(X \geq M) \leq \frac{\mathbb{E}X}{M} \leq \frac{cn}{M^{3/2}}$$

and setting $M = An^{2/3}$ for some large $A > 0$ concludes the proof.

The principle behind the proofs of this paper is that the tension between:

(a) $|\mathcal{C}(v)|$ is stochastically bounded by a critical branching process and
(b) $|\mathcal{C}(v)| \geq \beta n^{2/3}$

forces the geometry of $\mathcal{C}(v)$ to resemble that of a critical Galton–Watson tree conditioned to have at least $\beta n^{2/3}$ vertices. The following heuristic argument guides many of the proofs in this paper. Let $\mathcal{A}$ be a property of graphs (such as large volume or large diameter) which is inherited from a component $\mathcal{C}(v)$ by its bounding tree $T$. An argument in the same spirit as the above then gives that for $M = \beta n^{2/3}$, we have

$$\mathbf{P}(T \in \mathcal{A} \mid |T| \geq M) = o(1) \implies \mathbf{P}(|\mathcal{C}_1| \geq M \text{ and } \mathcal{C}_1 \in \mathcal{A}) = o(1).$$

Rigorous instances of this heuristic can be found in the proofs of all the theorems of this paper.

**2. Preliminaries.** Given a graph $G = (V, \mathcal{E})$ and a vertex $v$, denote by $\mathcal{C}(v) = \mathcal{C}(v, G_p)$ the connected component of $G_p$ which contains $v$. For a set of vertices $V' \subset V$, we write $\mathcal{E}(V')$ for the set of edges which have both ends in $V'$. We write $d_p(v, u)$ for the graph distance between $v$ and $u$ in $G_p$ and we define

$$B_p(v, k) = \{u \in \mathcal{C}(v) : d_p(v, u) \leq k\},$$
$$\partial B_p(v, k) = \{u : d_p(v, u) = k\}.$$

For ease of exposition, we begin by proving Theorem 1.2 assuming $\lambda \leq 0$. The case $\lambda > 0$ is proved in Section 6. Theorem 1.2 with $\lambda \leq 0$ will follow from the following, more general, theorem.

THEOREM 2.1. *Let $G = (V, \mathcal{E})$ be a graph and let $p \in (0, 1)$. Suppose that for some constants $c_1, c_2 > 0$ and all vertices $v \in V$, the following two conditions are satisfied for any subgraph $G' \subset G$:*

(i) $\mathbb{E}|\mathcal{E}(B_p(v, k))| \leq c_1 k$;
(ii) $\mathbf{P}(|\partial B_p(v, k)| > 0) \leq c_2/k$.

*Then:*

(a) $\mathbf{P}(\exists \mathcal{C} \in \mathbf{CO}(G_p) \text{ with } |\mathcal{C}| > \beta n^{2/3} \text{ and } \operatorname{diam}(\mathcal{C}) \notin [A^{-1} n^{1/3}, A n^{1/3}]) \leq O(A^{-1})$,
(b) $\mathbf{P}(\exists \mathcal{C} \in \mathbf{CO}(G_p) \text{ with } |\mathcal{E}(\mathcal{C})| > A n^{2/3}) \leq O(A^{-1})$,
(c.1) $\mathbf{P}(\exists \mathcal{C} \in \mathbf{CO}(G_p) \text{ with } |\mathcal{C}| > \beta n^{2/3} \text{ and } T_{\mathrm{mix}}(\mathcal{C}) > An) \leq O(A^{-1/2})$,
(c.2) $\mathbf{P}(\exists \mathcal{C} \in \mathbf{CO}(G_p) \text{ with } |\mathcal{C}| > \beta n^{2/3} \text{ and } T_{\mathrm{mix}}(\mathcal{C}) < A^{-1} n) \leq O(A^{-1/13})$,

*where the constants implicit in the O-notation depend on $c_1$, $c_2$ and $\beta$.*

PROOF OF THEOREM 1.2 FOR $\lambda \leq 0$. We verify the assumptions of Theorem 2.1 for a graph $G$ with maximum degree at most $d$ and $p \leq \frac{1}{d-1}$; we then take $A$ large enough. This is done by bounding the breadth–first search in the component of a vertex $v$ in $G_p$ by a breadth first search in a random tree. Let $\Gamma$ be an infinite $d$-regular tree with root $\rho$ (i.e., $\rho$ has $d$ children in the tree and any other vertex has $d - 1$ children and one parent) and let $d_\Gamma(u, v)$ denote the distance between vertices $u$ and $v$ in $\Gamma$. We denote by $\mathcal{C}(\rho,\Gamma_p)$ the component of $\rho$ in the subgraph $\Gamma_p$ obtained from percolation on $\Gamma$; let $\mathcal{L}_k$ be the set of vertices in level $k$ of $\mathcal{C}(\rho,\Gamma_p)$, that is,

$$\mathcal{L}_k = \{u \in \mathcal{C}(\rho,\Gamma_p) : d_\Gamma(\rho,u) = k\}.$$

Since the maximal degree in $G$ is at most $d$, we can clearly couple $G_p$ and $\Gamma_p$ so that the following two conditions hold:

(1) $|B_p(v,k)| \leq |\mathcal{E}(B_p(v,k))| + 1 \leq \sum_{j=0}^{k} |\mathcal{L}_j|$;
(2) $|\partial B_p(v,k)| \leq |\mathcal{L}_k|$.

Since

$$\mathbb{E}|\mathcal{L}_k| = d(d-1)^{k-1} p^k \leq 2$$

for all $k$, condition (i) of Theorem 2.1 is satisfied with $c_1 = 2$. For condition (ii), we use the following result, due to Lyons [19], which is related to an asymptotic estimate of Kolmogorov [16] (see also [15] and [20] for refinements and alternative proofs).

LEMMA 2.2 (Theorem 2.1 of [19]). *Assign each edge $e$ from level $k-1$ to level $k$ of $\Gamma$, the edge resistance $r_e = \frac{1-p}{p^k}$. Let $\mathcal{R}_k$ be the effective resistance from the root to level $k$ of $\Gamma$. Then,*

$$(2.1) \qquad \mathbf{P}(\mathcal{L}_k \neq \varnothing) \leq \frac{2}{1 + \mathcal{R}_k}.$$

Since $p \leq \frac{1}{d-1}$ and the edge resistances $r_e$ are monotone decreasing in $p$, the effective resistance $\mathcal{R}_k$ from $\rho$ to level $k$ of $\Gamma$ satisfies (see [24], Example 8.3)

$$\mathcal{R}_k = \sum_{i=1}^{k} \frac{(1-p)p^{-i}}{d(d-1)^{i-1}} \geq \sum_{i=1}^{k} \frac{\frac{d-2}{d-1}(d-1)^i}{d(d-1)^{i-1}} \geq \frac{(d-2)k}{d} \geq \frac{k}{3}$$

as $d \geq 3$. Thus, by our coupling and Lemma 2.2, condition (ii) holds with $c_2 = 6$. $\square$

### 3. The diameter of critical random graphs.

PROOF OF THEOREM 2.1(a). If a vertex $v \in V$ satisfies $\text{diam}(\mathcal{C}(v)) > R$, then $|\partial B_p(v, \lceil R/2 \rceil)| > 0$, hence, by condition (ii), we have

$$(3.1) \qquad \mathbf{P}(\text{diam}(\mathcal{C}(v)) > R) \leq \frac{2c_2}{R}.$$

If we write

$$X = |\{v \in V : \text{diam}(\mathcal{C}(v)) > R\}|,$$

then (3.1) implies that $\mathbb{E}X \leq \frac{2c_2 n}{R}$. By Markov's inequality, we have

$$\mathbf{P}(\exists \mathcal{C} \in \mathbf{CO}(G_p) \text{ with } |\mathcal{C}| > M \text{ and } \operatorname{diam}(\mathcal{C}) > R)$$
(3.2)
$$\leq \mathbf{P}(X > M) \leq \frac{2c_2 n}{MR}.$$

If $v \in V$ satisfies $\operatorname{diam}(\mathcal{C}(v)) \leq r$ and $|\mathcal{C}(v)| > M$, then $|\mathcal{E}(B_p(v,r))| \geq M$. Thus, by condition (i) and Markov's inequality, we have

(3.3) $$\mathbf{P}(\operatorname{diam}(\mathcal{C}(v)) \leq r \text{ and } |\mathcal{C}(v)| > M) \leq \frac{c_1 r}{M}.$$

If we write

$$Y = |\{v \in V : |\mathcal{C}(v)| > M \text{ and } \operatorname{diam}(\mathcal{C}(v)) < r\}|,$$

then (3.3) implies that $\mathbb{E}Y \leq \frac{c_1 r n}{M}$, whence, by Markov's inequality,

$$\mathbf{P}(\exists \mathcal{C} \in \mathbf{CO}(G_p) \text{ with } |\mathcal{C}| > M \text{ and } \operatorname{diam}(\mathcal{C}) < r)$$
(3.4)
$$\leq \mathbf{P}(Y > M) \leq \frac{c_1 r n}{M^2}.$$

Combining (3.2) and (3.4) gives

$$\mathbf{P}(\exists \mathcal{C} \in \mathbf{CO}(G_p) \text{ with } |\mathcal{C}| > M \text{ and } \operatorname{diam}(\mathcal{C}) \notin [r, R]) \leq \left(\frac{c_1 r}{M} + \frac{2c_2}{R}\right)\frac{n}{M}.$$

Take $M = \beta n^{2/3}$ and set $r = A^{-1} n^{1/3}$ and $R = A n^{1/3}$. The right-hand side of the preceding display is then $(c_1 \beta^{-2} + 2c_2 \beta^{-1}) A^{-1} = O(A^{-1})$, which completes the proof. □

PROOF OF THEOREM 2.1(b). In this proof, we will only use conditions (i) and (ii) of the theorem for $k \leq n^{1/3}$. Fix some $M > 1$ and $r \leq n^{1/3}$. Observe that for any $v \in V$, we have

$$\{|\mathcal{C}(v)| > M\} \subset \{|\mathcal{C}(v)| > M \text{ and } \operatorname{diam}(\mathcal{C}(v)) \leq r\} \cup \{\operatorname{diam}(\mathcal{C}(v)) > r\}.$$

Write

$$X = |\{v \in V : |\mathcal{C}(v)| > M\}|.$$

By condition (ii), we have that (3.1) holds for $R < 2n^{1/3}$ and by condition (i), we have that (3.3) holds for $r < n^{1/3}$. Thus, by taking $R = r$ in (3.1) and (3.3), we deduce that

$$\mathbb{E}X \leq \left(\frac{2c_2}{r} + \frac{c_1 r}{M}\right) n.$$

Let $\mathcal{C}_1$ denote the largest component of $G_p$. Observe that if $|\mathcal{C}_1| \geq M$, then $|\mathcal{C}_1| \leq X$. Therefore, $|\mathcal{C}_1| \leq M + X$. We take $M = \lceil n^{2/3} \rceil$ and $r = \lfloor n^{1/3} \rfloor$ and obtain that $\mathbb{E}|\mathcal{C}_1| \leq (2c_1 + c_2 + 2)n^{2/3}$. We then have that for any $\widetilde{A} > 0$,

$$\mathbf{P}(\exists \mathcal{C} \in \mathbf{CO}(G_p) \text{ with } |\mathcal{C}| \geq \widetilde{A}n^{2/3}) \leq O(\widetilde{A}^{-1}). \tag{3.5}$$

Next, observe that condition (i) for $k = 1$ implies that the maximal degree $d$ in $G$ satisfies $dp \leq c_1$. Consider "exploring" the levels $\partial B_p(v, k)$ level by level. At the end, we have discovered a spanning tree on the vertices of $\mathcal{C}(v)$ and since $d$ is the maximal degree, the number of extra edges in this component can be bounded above by $Z$, a random variable distributed as $\mathrm{Bin}(d|\mathcal{C}(v)|, p)$. Thus, if we condition on the vertices of $\mathcal{C}(v)$, then the number of edges $|\mathcal{E}(\mathcal{C}(v))|$ can be stochastically bounded above by $|\mathcal{C}(v)| - 1 + Z$. By a standard large deviation inequality (see, e.g., [3], Section A.14), we have $\mathbf{P}(Z \geq 2c_1 m | |\mathcal{C}(v)| = m) \leq e^{-2c_1 \gamma m}$ for some constant $\gamma > 0$, so

$$\mathbf{P}\left(|\mathcal{E}(\mathcal{C}(v))| \geq An^{2/3} \text{ and } |\mathcal{C}(v)| < \frac{A}{2c_1 + 1}n^{2/3}\right) \leq e^{-\gamma An^{2/3}}.$$

Thus,

$$\mathbf{P}\left(\exists \mathcal{C} \in \mathbf{CO}(G_p) \text{ with } |\mathcal{E}(\mathcal{C})| \geq An^{2/3} \text{ and } |\mathcal{C}| \leq \frac{A}{2c_1 + 1}n^{2/3}\right) \leq ne^{-\gamma An^{2/3}}.$$

This, together with (3.5) gives that

$$\mathbf{P}(\exists \mathcal{C} \in \mathbf{CO}(G_p) \text{ with } |\mathcal{E}(\mathcal{C})| \geq An^{2/3}) \leq O(A^{-1}),$$

concluding our proof. □

**4. The upper bound on the mixing time.** The following known lemma bounds the total variation mixing time in terms of the maximal hitting time. For variants of this lemma, see Chapter 4 of [2].

LEMMA 4.1. *Let $\mathbf{p}$ be transition probabilities for a reversible, lazy [i.e., $\mathbf{p}(x, x) \geq 1/2$ for all $x \in V$] Markov chain on a finite state space $V$. For $x \in V$, denote by $\tau_x$ the hitting time of $x$. We have*

$$T_{\mathrm{mix}}(1/4) \leq 2 \max_{x,y \in V} \mathbb{E}_y \tau_x.$$

PROOF. Lemma 11 of Chapter 2 in [2] states that

$$\pi(x)\mathbb{E}_\pi(\tau_x) = \sum_{t=0}^{\infty} [\mathbf{p}^t(x, x) - \pi(x)],$$

where $\pi$ is the stationary distribution. Let $\{\lambda_i\}_{i=1}^{|V|}$ be the eigenvalues of the transition matrix $\mathbf{p}$, with corresponding (real) eigenfunctions $\{\psi_i\}_{i=1}^{|V|}$,

normalized in $L^2(\pi)$. (In particular, $\lambda_1 = 1$ and $\psi_1 \equiv 1$.) By spectral decomposition, for each $x \in V$,

$$\mathbf{p}^t(x,x) = \pi(x) \sum_{i=1}^{|V|} \lambda_i^t \psi_i(x)^2.$$

Since the chain is lazy, $\lambda_i \in [0,1]$ for all $i$ and hence $\mathbf{p}^{t+1}(x,x) \leq \mathbf{p}^t(x,x)$ for all $t \geq 0$. Therefore, for any integer $m > 0$,

$$\pi(x) \mathbb{E}_\pi(\tau_x) \geq \sum_{t=1}^{2m} [\mathbf{p}^t(x,x) - \pi(x)] \geq 2m[\mathbf{p}^{2m}(x,x) - \pi(x)],$$

hence

(4.1) $$\frac{\mathbb{E}_\pi(\tau_x)}{2m} \geq \frac{\mathbf{p}^{2m}(x,x)}{\pi(x)} - 1.$$

By Cauchy–Schwarz, we have

$$4\|\mathbf{p}^m(x,\cdot) - \pi\|_{\mathrm{TV}}^2 = \left( \sum_y \pi(y) \left| \frac{\mathbf{p}^m(x,y)}{\pi(y)} - 1 \right| \right)^2$$

$$\leq \sum_y \pi(y) \left[ \frac{\mathbf{p}^m(x,y)}{\pi(y)} - 1 \right]^2,$$

therefore, by reversibility, we obtain

$$4\|\mathbf{p}^m(x,\cdot) - \pi\|_{\mathrm{TV}}^2 \leq \sum_y \left[ \frac{\mathbf{p}^m(x,y)\mathbf{p}^m(y,x)}{\pi(x)} - 2\mathbf{p}^m(x,y) + \pi(y) \right]$$

$$= \frac{\mathbf{p}^{2m}(x,x)}{\pi(x)} - 1.$$

Thus, by (4.1), we obtain

$$4\|\mathbf{p}^m(x,\cdot) - \pi\|_{\mathrm{TV}}^2 \leq \frac{\mathbb{E}_\pi(\tau_x)}{2m} \leq \frac{\max_{x,y \in V} \mathbb{E}_y(\tau_x)}{2m}$$

and the right-hand side is at most $\frac{1}{4}$ for $m \geq 2 \max_{x,y \in V} \mathbb{E}_y(\tau_x)$, concluding our proof. □

REMARK. Lemma 4.1 actually gives a bound on the $\ell^2$-mixing time.

COROLLARY 4.2. *Let $G = (V, \mathcal{E})$ be a graph. The mixing time of a lazy simple random walk on $G$ then satisfies*

$$T_{\mathrm{mix}}(G, 1/4) \leq 8|\mathcal{E}(G)|\mathrm{diam}(G).$$

PROOF. For any two vertices $x$ and $y$, let $d_G(x,y)$ denote the graph distance in $G$ between $x$ and $y$. We bound $\mathbb{E}_y(\tau_x)$ by $\mathbb{E}_y(\tau_x) + \mathbb{E}_x(\tau_y)$, which is also known as the *commute time* between $x$ and $y$. Let $\mathcal{R}(x \leftrightarrow y)$ denote the effective resistance from $x$ to $y$ when each edge has unit resistance. The commute time identity of [8] (see also [25]) implies that for lazy simple random walk on a connected graph $G = (V, \mathcal{E})$,

$$\mathbb{E}_y(\tau_x) + \mathbb{E}_x(\tau_y) = 4|\mathcal{E}(G)|\mathcal{R}(x \leftrightarrow y).$$

Since $\mathcal{R}(x \leftrightarrow y) \leq d_G(x,y)$, Lemma 4.1 concludes the proof. □

We are now ready to prove the *mixing time upper bound.*

PROOF OF THEOREM 2.1(c.1). If a cluster $\mathcal{C} \in \mathbf{CO}(G_p)$ satisfies $T_{\mathrm{mix}}(\mathcal{C}) > An$, then $|\mathcal{E}(\mathcal{C})|\mathrm{diam}(\mathcal{C}) > (A/8)n$, by Corollary 4.2, so either $|\mathcal{E}(\mathcal{C})| > (A/8)^{1/2}n^{2/3}$ or $\mathrm{diam}(\mathcal{C}) > (A/8)^{1/2}n^{1/3}$. By part (b) of Theorem 2.1, we have

$$\mathbf{P}(\exists \mathcal{C} \in \mathbf{CO}(G_p) : |\mathcal{E}(\mathcal{C})| \geq (A/8)^{1/2}n^{2/3}) \leq O(A^{-1/2})$$

and by part (a) of Theorem 2.1 we have

$$\mathbf{P}(\exists \mathcal{C} \in \mathbf{CO}(G_p) : \mathrm{diam}(\mathcal{C}) > (A/8)^{1/2}n^{1/3}) \leq O(A^{-1/2}).$$

Adding the probabilities in the last two displays proves the proposition. □

**5. The lower bound on the mixing time.** We will use the Nash-Williams inequality [23] (see also [24]). Recall that a set of edges $\Pi$ is a *cut-set* separating a vertex $v$ from a set of vertices $U$ if any path from $v$ to $U$ must intersect $\Pi$.

LEMMA 5.1 (Nash-Williams [23]). *If $\{\Pi_j\}_{j=1}^J$ are disjoint cut-sets separating $v$ from $U$ in a graph with unit conductance for each edge, then the effective resistance from $v$ to $U$ satisfies*

$$\mathcal{R}(v \leftrightarrow U) \geq \sum_{j=1}^{J} \frac{1}{|\Pi_j|}.$$

We will also use the following lemma, due to Tetali [25].

LEMMA 5.2 (Tetali [25]). *Let $\mathbf{p}$ be transition probabilities for a reversible Markov chain. Let $\mu : V \to \mathbb{R}$ be a function such that $\mu(x)\mathbf{p}(x,y) = \mu(y)\mathbf{p}(y,x)$ for all $x, y$, and let $c_{x,y} = \mu(x)\mathbf{p}(x,y)$ be the edge conductances. Then, for this Markov chain,*

$$\mathbb{E}_v(\tau_z) = \tfrac{1}{2} \sum_{u \in V} \mu(u)[\mathcal{R}(v \leftrightarrow z) + \mathcal{R}(z \leftrightarrow u) - \mathcal{R}(u \leftrightarrow v)].$$

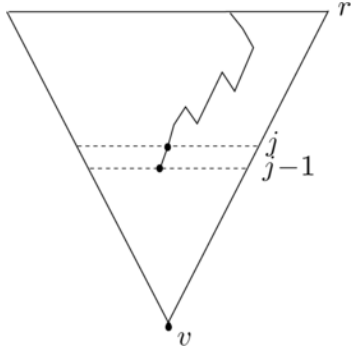

Fig. 1. *A lane.*

COROLLARY 5.3. *For a lazy random walk on a finite graph where each edge has unit resistance, we have*

$$\mathbb{E}_v(\tau_z) = \sum_{u \in V} \deg(u)[\mathcal{R}(v \leftrightarrow z) + \mathcal{R}(z \leftrightarrow u) - \mathcal{R}(u \leftrightarrow v)].$$

PROOF. Take $\mu(x) = 2\deg(x)$ so that $\mu(x)\mathbf{p}(x,y) = 1$. □

The following structural argument is inspired by Barlow and Kumagai [4]. For a graph $G = (V, \mathcal{E})$, write $d_G(x, y)$ for the graph distance between $x$ and $y$. For any vertex $v$, let

$$B(v, r) = B_G(v, r) = \{u \in v : d_G(u, v) \leq r\},$$
$$\partial B(v, r) = \partial B_G(v, r) = \{u \in v : d_G(u, v) = r\}.$$

To motivate the following *definitions*, think of the edges of $B(v, r)$ as a road network that connects $v$ to $\partial B(v, r)$.

- An edge $e$ between $\partial B(v, j-1)$ and $\partial B(v, j)$ is called a *lane* for $(v, r)$ if it there is a path with initial edge $e$ from $\partial B(v, j-1)$ to $\partial B(v, r)$ that does not return to $\partial B(v, j-1)$.
- we say that a level $j$ (with $0 < j < r$) has $L$ lanes for $(v, r)$ if there are at least $L$ edges between $\partial B(v, j-1)$ and $\partial B(v, j)$ which are lanes for $(v, r)$.
- Let $k < r$. A vertex $v$ is called *L-lane rich* for $(k, r)$ if more than half of the levels $j \in [k/2, k]$ have $L$ lanes for $(v, r)$.

LEMMA 5.4. *Let $G = (V, \mathcal{E})$ be a graph and let $v \in V$. Suppose that $|B(v, h)| \geq m$, that $v$ is not L-lane rich for $(k, r)$, that $|\mathcal{E}(B(v, r))| < \frac{|\mathcal{E}(G)|}{3}$ and that $h < \frac{k}{4L}$. Then,*

$$T_{\mathrm{mix}}(G) \geq \frac{mk}{12L}.$$

PROOF. As $v$ is not $L$-lane rich for $(k,r)$, there are at least $k/4$ levels between $k/2$ and $k$ which have less than $L$ lanes for $(v,r)$. In each such level $j$, the lanes for $(v,r)$ form a cut-set of size less than $L$ separating any $u \in B(v,h)$ from $\partial B(v,r)$. Thus, for any $u \in B(v,h)$, the Nash-Williams inequality, Lemma 5.1, yields

$$\mathcal{R}(u \leftrightarrow \partial B(v,r)) \geq \frac{k}{4L}. \tag{5.1}$$

By the triangle inequality for effective resistance (see, e.g., [25]), each of the summands of Lemma 5.3 is nonnegative. Denote by $\tau[r]$ the hitting time of $\partial B(v,r)$ by the lazy simple random walk. By gluing $\partial B(v,r)$ into a single vertex $z$, Corollary 5.3 gives

$$\mathbb{E}_v \tau[r] \geq \sum_{u \in B(v,h)} [\mathcal{R}(v \leftrightarrow \partial B(v,r)) + \mathcal{R}(u \leftrightarrow \partial B(v,r)) - \mathcal{R}(u \leftrightarrow v)];$$

by (5.1) and the fact that $\mathcal{R}(u \leftrightarrow v) \leq d_G(u,v)$, we infer that

$$\mathbb{E}_v \tau[r] \geq |B(v,h)| \left( \frac{k}{2L} - h \right) \geq \frac{k}{4L} |B(v,h)|,$$

where the last inequality is due to our assumption on $h$. Since $|B(v,h)| \geq m$,

$$\mathbb{E}_v \tau[r] \geq \frac{mk}{4L}. \tag{5.2}$$

Fix some integer $t > 0$. If all vertices $x \in B(v, r-1)$ satisfy $\mathbf{P}_x(\tau[r] \leq t) \geq \frac{1}{3}$, then $\tau[r]/t$ is stochastically dominated by a geometric$(1/3)$ random variable, whence $\mathbb{E}_v(\tau[r]) \leq 3t$. By this and (5.2), we conclude that for $t = \frac{mk}{12L}$, there exists some $x \in B(v, r-1)$ such that $\mathbf{P}_x(\tau[r] \leq t) \leq \frac{1}{3}$. Therefore, for this $t$, we have $\mathbf{p}^t(x, B(v,r)) \geq \frac{2}{3}$ and as $|\mathcal{E}(B(v,r))| < |\mathcal{E}(G)|/3$, we have $\pi(B(v,r)) \leq 1/3$. We thus have that $\|\mathbf{p}^t(x,\cdot) - \pi(\cdot)\|_{\mathrm{TV}} > \frac{1}{4}$ and therefore

$$T_{\mathrm{mix}}(G) \geq \frac{mk}{12L}. \qquad \square$$

We return to the setting of Theorem 2.1 and write $B_p(v,r)$ for $B_{\mathcal{C}(v)}(v,r)$. Define:

- $\mathcal{A}_1(v,h,m) = \{|B_p(v,h)| < m\}$;
- $\mathcal{A}_2(v,L,k,r) = \{v \text{ is } L\text{-lane rich for } (k,r) \text{ in } \mathcal{C}(v)\}$;
- $\mathcal{A}_3(v,\alpha,r) = \{|\mathcal{E}(B_p(v,r))| \geq \alpha r^2\}$.

PROPOSITION 5.5. *Suppose that conditions* (i) *and* (ii) *of Theorem* 2.1 *hold for all* $k \leq n^{1/3}$. *If* $h < \frac{n^{1/3}}{4}$, *then*

$$\mathbf{P}(|\mathcal{C}(v)| > \beta n^{2/3} \text{ and } \mathcal{A}_1(v,h,m)) \leq \frac{4mc_2^2}{h^3} + \frac{4c_1 h}{\beta n^{2/3}}.$$

PROOF. If $|B_p(v,h)| < m$, then there exists a level $j \in [h/2, h]$ such that $|\partial B_p(v,j)| \leq \frac{2m}{h}$. Fix the smallest such $j$. If, in addition, $\mathrm{diam}(\mathcal{C}(v)) > 4h$, then $\partial B_p(v, 2h) \neq \varnothing$, so at least one of the at most $\frac{2m}{h}$ vertices $w$ in $\partial B_p(v,j)$ must be the beginning of a path in $\mathcal{C}(v)$ that does not return to level $j$ and reaches at least $2h - j \geq h$ levels higher; given $w$, the existence of such a path has probability at most $c_2/h$ by condition (ii). Applying (ii) again, together with the Markov property at level $j$, we deduce that

$$\mathbf{P}(\mathrm{diam}(\mathcal{C}(v)) > 4h \text{ and } |B_p(v,h)| < m) \leq \mathbf{P}\left(\partial B_p\left(v, \frac{h}{2}\right) \neq \varnothing\right) \cdot \frac{2m}{h} \cdot \frac{c_2}{h}$$

$$\leq \frac{c_2}{h/2} \cdot \frac{2m}{h} \cdot \frac{c_2}{h} = \frac{4mc_2^2}{h^3}.$$

Since (3.3) holds for $r < n^{1/3}$, by condition (i) and $4h < n^{1/3}$, combining it with the last display gives

$$\mathbf{P}(|\mathcal{C}(v)| > \beta n^{2/3} \text{ and } |B_p(v,h)| < m) \leq \frac{4mc_2^2}{h^3} + \frac{4c_1 h}{\beta n^{2/3}}. \qquad \square$$

PROPOSITION 5.6. *Under the conditions of Theorem 2.1, if $k \leq r/2$ and $r < n^{1/3}$, then*

$$\mathbf{P}(\mathcal{A}_2(v, L, k, r)) \leq \frac{8c_1 c_2}{Lr}.$$

PROOF. For each edge between $\partial B_p(v, j-1)$ and $\partial B_p(v, j)$, where $j \leq k$, the probability that it begins a path to $\partial B_p(v, r)$ that does not go through $\partial B_p(v, j-1)$ is at most $\frac{c_2}{r-j}$, by condition (ii). This, with condition (i), implies that the expected number of lanes for $(v, r)$ in $\mathcal{E}(B_p(v, k))$ is at most $\frac{c_1 c_2 k}{r-j}$. If $v$ is $L$-lane rich for $(k, r)$, then there are at least $\frac{Lk}{4}$ lanes in $\mathcal{E}(B_p(v, k))$. Thus, Markov's inequality gives

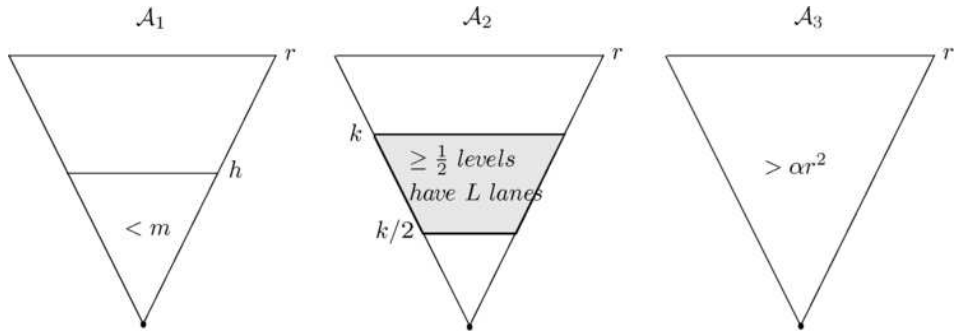

FIG. 2.

$$\mathbf{P}(\mathcal{A}_2(v,L,k,r)) \leq \frac{4c_1c_2k}{(r-j)Lk} \leq \frac{8c_1c_2}{Lr},$$

by our assumption on $k$. □

PROPOSITION 5.7. *Under the conditions of Theorem* 2.1, *if* $r < n^{1/3}$, *then*

$$\mathbf{P}(\mathcal{A}_3(v,\alpha,r)) \leq \frac{c_1}{\alpha r}.$$

PROOF. This follows by Markov's inequality from condition (i). □

We are now ready to prove the *mixing time lower bound*.

PROOF OF THEOREM 2.1(c.2). We abbreviate $\mathcal{A}_1(v,h,m)$, $\mathcal{A}_2(v,L,k,r)$ and $\mathcal{A}_3(v,\alpha,r)$ by $\mathcal{A}_1, \mathcal{A}_2$ and $\mathcal{A}_3$, respectively. Fix $D > 0$, to be chosen later, and define the following parameters:

$$L = \beta^{-3}D^2, \qquad \alpha = \frac{\beta^{-3}D^2}{20}, \qquad h = \frac{\beta^5 D^{-3}}{4} n^{1/3},$$

$$k = 5Lh, \qquad r = 10Lh, \qquad m = h^3\beta D^{-1}n^{-1/3}.$$

By Propositions 5.5, 5.6 and 5.7, we have that

$$\mathbf{P}(|\mathcal{C}(v)| > \beta n^{2/3} \text{ and } (\mathcal{A}_1 \cup \mathcal{A}_2 \cup \mathcal{A}_3))$$

(5.3)
$$\leq \frac{4mc_2^2}{h^3} + \frac{4c_1h}{\beta n^{2/3}} + \frac{8c_1c_2}{Lr} + \frac{c_1}{\alpha r}$$

$$\leq (4c_2^2 + c_1\beta^3 D^{-2} + 4c_1c_2 + 8c_1)\beta D^{-1}n^{-1/3}.$$

Let

$$X = |\{v : |\mathcal{C}(v)| > \beta n^{2/3} \text{ and } (\mathcal{A}_1 \cup \mathcal{A}_2 \cup \mathcal{A}_3)\}|.$$

Then, by (5.3), we have $\mathbb{E}X \leq (4c_2^2 + c_1\beta^3 D^{-2} + 4c_1c_2 + 8c_1)\beta D^{-1}n^{2/3}$. Denote by $\mathcal{A}$ the event that there exists a component $\mathcal{C} \in \mathbf{CO}(G_p)$, such that $|\mathcal{C}| > \beta n^{2/3}$ and all of the vertices $v \in \mathcal{C}$ satisfy either $\mathcal{A}_1, \mathcal{A}_2$ or $\mathcal{A}_3$. Observe that if $\mathcal{A}$ holds, then $X > \beta n^{2/3}$. Thus, by Markov's inequality,

(5.4) $$\mathbf{P}(\mathcal{A}) \leq (4c_2^2 + c_1\beta^3 D^{-2} + 4c_1c_2 + 8c_1)D^{-1}.$$

If $\mathcal{A}$ does not hold, then all components $\mathcal{C}$ with $|\mathcal{C}| > \beta n^{2/3}$ have a vertex $v \in \mathcal{C}$ such that $|B_p(v,h)| \geq m$, the vertex $v$ is not $L$-lane rich for $(k,r)$ and $|\mathcal{E}(B_p(v,r))| < \alpha r^2$. It is easy to verify that $\beta n^{2/3} \geq \alpha r^2/3$ and that $h < \frac{k}{4L}$. Thus, Lemma 5.4 gives that, with probability at least $1 - \mathbf{P}(\mathcal{A})$, all components $\mathcal{C} \in \mathbf{CO}(G_p)$ with $|\mathcal{C}| > \beta n^{2/3}$ satisfy

$$T_{\mathrm{mix}}(\mathcal{C}) \geq \frac{mk}{12L} \geq \frac{\beta^{21}}{1000 D^{13}} n.$$

Setting $D = (\frac{A\beta^{21}}{1000})^{1/13}$, so that $A = 1000 D^{13} \beta^{-21}$, concludes the proof. □

**6. The diameter inside the scaling window.** The following theorem is essentially Theorem 2.1 under weaker conditions that hold for all $p \leq \frac{1+\lambda n^{-1/3}}{d-1}$ when $\lambda \in \mathbb{R}$ is fixed.

THEOREM 6.1. *Let $G = (V, \mathcal{E})$ be a graph and $p \in [0,1]$. Assume that the following holds for any subgraph $G' \subset G$:*

(i′) $\mathbb{E}|\mathcal{E}(B_p(v,k))| \leq c_1 k$ *for* $k \leq n^{1/3}$;
(ii′) $\mathbf{P}(|\partial B_p(v,k)| > 0) \leq c_2/k$ *for* $k \leq n^{1/3}$.

*The conclusions of Theorem 2.1 then hold and conclusions (1.1) and (1.2) of Theorem 1.3 also hold.*

PROOF OF THEOREM 1.2 FOR $\lambda \in \mathbb{R}$ AND THEOREM 1.3. We verify the assumptions of Theorem 6.1 by bounding the breadth-first search in the component of a vertex $v$ in $G_p$ by a breadth-first search in a random tree, as we did in the case $\lambda \leq 0$. We have

$$\mathbb{E}|\mathcal{L}_k| = d(d-1)^{k-1} p^k \leq 2(1 + \lambda n^{-1/3})^k \leq 2e^\lambda$$

for $k \leq n^{1/3}$ and thus, by the coupling from before, condition (i′) holds with $c_1 = 2e^\lambda$. In the notation of Lemma 2.2, we have

$$\mathcal{R}_k = \sum_{i=1}^k \frac{(1-p)p^{-i}}{d(d-1)^{i-1}}$$

$$\geq \sum_{i=1}^k \frac{\frac{d-2-\lambda n^{-1/3}}{d-1}(d-1)^i}{d(d-1)^{i-1}(1+\lambda n^{-1/3})^i}$$

$$\geq \left(\frac{d-2-\lambda n^{-1/3}}{d}\right) \frac{k}{(1+\lambda n^{-1/3})^k} \geq \frac{k}{4e^\lambda}$$

for $n$ large enough (as $d \geq 3$ and $k \leq n^{1/3}$). Thus, condition (ii) holds with $c_2 = 8e^\lambda$. □

The following lemmas will be essential for the proof of Theorem 6.1.

LEMMA 6.2. *Assume the conditions of Theorem* 6.1. *Let $M$ and $R$ be two positive integers satisfying*

$$R > 16c_2 M n^{-1/3}.$$

*Then,*

$$\mathbf{P}(|\mathcal{C}(v)| \leq M \text{ and } \mathrm{diam}(\mathcal{C}(v)) > 2R) \leq c_2\left(\frac{2}{R} \vee n^{-1/3}\right) 2^{-R^2/((64c_2+2)M)}.$$

PROOF. Set $h = M/R$. We may assume $2R < M$ as, otherwise, the required probability is 0. We say that level $j$ of the exploration tree from $v$ is *thin* if it contains at most $8h$ vertices. Define $j_1$ as the first thin level greater than $R/2$ and, for $i > 1$, define

(6.1) $$j_i := \min\{j > j_{i-1} + 16c_2 h : |\partial B_p(v,j)| \leq 8h\},$$

that is, $j_i$ is the first thin level greater than $j_{i-1} + 16c_2 h$.

We call a vertex $w \in \partial B_p(v,j)$ *good* if there is a path from $w$ to $\partial B_p(v, j + 16c_2 h)$ that intersects $B_p(v,j)$ only in $w$; we call level $j$ in the exploration process from $v$ *good* if it contains at least one good vertex. For each vertex $w \in \partial B_p(v,j)$, the conditional probability that it is good, given $B_p(v,j)$, is at most $\frac{1}{16h}$, by condition (ii′) (and the inequality $16c_2 h < n^{1/3}$, which follows from our assumption on $R$ and $M$). Therefore, for every $j$, we have that

$$\mathbf{P}(\text{level } j \text{ is good} \mid B_p(v,j), \text{level } j \text{ is thin}) \leq \tfrac{1}{2}.$$

By the previous display, we deduce that, with $j_i$ defined in (6.1), we have

(6.2) $$\mathbf{P}(\text{level } j_i \text{ is good for all } i \leq k-1 | B_p(v, \lceil R/2 \rceil)) \leq 2^{-(k-1)}.$$

If $|\mathcal{C}(v)| \leq M$ and $\mathrm{diam}(\mathcal{C}(v)) > 2R$, then levels $j_1, \ldots, j_{k-1}$ are good with

(6.3) $$k - 1 \geq \frac{R}{(64c_2 + 2)h}.$$

To see this, let $\ell$ be the number of thin levels $j$ such that $\frac{R}{2} \leq j \leq R$. Since $|C(v)| \leq M$, we must have that $\ell \geq \frac{R}{4}$. From these $\ell \geq R/4$ levels, we obtain, in (6.1), at least $\frac{R}{4(16c_2 h + 1)} > \frac{R}{(64c_2+2)h}$ (as $h \geq 2$) thin levels separated from each other by more than $16c_2 h$ levels.

If $|C(v)| \leq M$ and $\mathrm{diam}(\mathcal{C}(v)) > 2R$, then $|\{\partial B_p(v, \lceil R/2 \rceil)\}| > 0$. By condition (ii′) of Theorem 6.1, we have

(6.4) $$\mathbf{P}(|\partial B_p(v, \lceil R/2 \rceil)| > 0) \leq c_2\left(\frac{2}{R} \vee n^{-1/3}\right).$$

Thus, if $|\mathcal{C}(v)| \leq M$ and $\mathrm{diam}(\mathcal{C}(v)) > 2R$, then the event in the left-hand side of (6.2) also holds, where $k$ satisfies (6.3). Therefore,

$$\mathbf{P}(|\mathcal{C}(v)| \leq M \text{ and } \mathrm{diam}(\mathcal{C}(v)) > 2R)$$
$$\leq c_2 \left(\frac{2}{R} \vee n^{-1/3}\right) 2^{-R^2/((64c_2+2)M)}. \qquad \square$$

LEMMA 6.3. *Assume the conditions of Theorem* 6.1. *Let $M$ and $R$ be two positive integers satisfying $R > 32c_2 M n^{-1/3}$ and $R > \sqrt{4(64c_2+2)M}$. Then,*

$$\mathbf{P}(\exists \mathcal{C} \in \mathbf{CO}(G_p) \text{ with } |\mathcal{C}| \leq M \text{ and } \mathrm{diam}(\mathcal{C}) > 2R)$$
$$\leq 4c_2 \left(\frac{2}{R} \vee n^{-1/3}\right) \frac{2^{-R^2/(2(64c_2+2)M)}}{M} n.$$

PROOF. Let $X_m$ be the random variable

$$X_m = \left|\left\{v \in V : \frac{m}{2} \leq |\mathcal{C}(v)| \leq m \text{ and } \mathrm{diam}(\mathcal{C}(v)) > 2R\right\}\right|.$$

By Lemma 6.2 and our assumptions on $M$ and $R$, for every $1 \leq k \leq \lceil \log_2(M) \rceil$, we have

$$\mathbb{E} X_{2^k} \leq n c_2 \left(\frac{2}{R} \vee n^{-1/3}\right) 2^{-R^2/((64c_2+2)2^k)}.$$

Let $\mathcal{A}$ denote the event

$$\mathcal{A} = \{\exists \mathcal{C} \in \mathbf{CO}(G_p) \text{ with } |\mathcal{C}| \leq M \text{ and } \mathrm{diam}(\mathcal{C}) > 2R\}.$$

$\mathcal{A}$ then implies that $X_{2^k} \geq 2^{k-1}$ for at least one $k$ satisfying $1 \leq k \leq \lceil \log_2(M) \rceil$. By applying Markov's inequality, we obtain

$$\mathbf{P}(\mathcal{A}) \leq n c_2 \left(\frac{2}{R} \vee n^{-1/3}\right) \sum_{k=1}^{\lceil \log_2 M \rceil} 2^{-R^2/((64c_2+2)2^k) - k + 1}.$$

It is straightforward to check that since $R > \sqrt{4(64c_2+2)M}$, the $k$th summand in the above sum is at most $1/2$ of the next summand, whence

$$\mathbf{P}(\mathcal{A}) \leq 4c_2 \left(\frac{2}{R} \vee n^{-1/3}\right) \frac{2^{-R^2/(2(64c_2+2)M)}}{M} n. \qquad \square$$

PROOF OF THEOREM 6.1. In the proofs of part (b) and part (c.2) of Theorem 2.1, we only used the weaker conditions (i') and (ii') of Theorem 6.1 [rather than (i) and (ii)], so no additional work is required there. Also, (3.4) holds for $r < n^{1/3}$, so taking $r = A^{-1} n^{1/3}$ gives the lower bound on the diameter implied in part (a) of the theorem. By Corollary 4.2, part (c.1) is an immediate corollary of the upper bound on the diameter (1.1) and part (b) of Theorem 2.1. Thus, all that is left to prove is (1.1) and (1.2).

PROOF OF (1.1). Take large $A$ and set $R = \lceil An^{1/3} \rceil$ and $M = \lfloor \frac{An^{2/3}}{32c_2} \rfloor$. Note that the assumptions of Lemma 6.3 are satisfied. Thus, part (b) of the theorem and Lemma 6.3 with these chosen $R$ and $M$ gives that

$$\mathbf{P}(\exists \mathcal{C} \in \mathbf{CO}(G_p) \text{ with } \operatorname{diam}(\mathcal{C}) > An^{1/3}) \leq O(A^{-1}),$$

which finishes the proof of (1.1).

PROOF OF (1.2). Take $M \leq n^{2/3}/2$ and $R = 2D\sqrt{M \log(n/M^{3/2})}$ for some large $D$ and substitute into Lemma 6.3 to prove (1.2). □

PROOF OF PROPOSITION 1.4. It is proved in [21] and [22] that under the assumptions of the proposition,

(6.5) $$\mathbf{P}(\exists \mathcal{C} \in \mathbf{CO}(G_p) \text{ with } |\mathcal{C}| > Bn^{2/3}) \leq e^{-\gamma B^3}$$

for some $\gamma > 0$ that depends on $\lambda$ (and on $d$ for the case of the random $d$-regular graph). For large enough $A$, take $R = \lfloor An^{1/3} \rfloor$ and $M = A^{1/2}n^{2/3}$; then, by Lemma 6.3,

$$\mathbf{P}(\exists \mathcal{C} \in \mathbf{CO}(G_p) \text{ with } |\mathcal{C}| \leq M \text{ and } \operatorname{diam}(\mathcal{C}) > 2R)$$

$$\leq 4c_2 n^{2/3} \frac{2^{-R^2/(2(64c_2+2)M)}}{M} \leq 2^{-\delta A^{3/2}},$$

where $\delta > 0$ is an absolute constant. In conjunction with (6.5) for $B = A^{1/2}$, this gives

$$\mathbf{P}(\exists \mathcal{C} \in \mathbf{CO}(G_p) \text{ with } \operatorname{diam}(\mathcal{C}) > An^{1/3}) \leq e^{-cA^{3/2}},$$

for some $c > 0$. □

## 7. Concluding remarks.

1. Theorems 2.1 and 6.1 naturally lead to the following question: *for which graphs $G$ and retention probabilities $p$ are conditions* (i'), (ii') *of Theorem 6.1 satisfied, yet there is a substantial probability of having connected components of size $n^{2/3}$?* In particular, it seems interesting to prove that these conditions hold for the Hamming cube $\{0,1\}^n$ or the $d$-dimensional discrete torus $[n]^d$ for large $d$ and some $p$.
2. A far more challenging problem is finding a good definition for the *critical probability* for percolation on finite transitive graphs. Among other properties, we expect that at this critical probability, the size of the largest component is not concentrated and that the second largest component has maximal expectation. We also expect to find a scaling window around this critical probability in which the above properties still hold, while outside this window, the properties cease to hold.

The following is a suggestion for such a definition. Let $\chi(p) = \mathbb{E}_p|\mathcal{C}(v)|$. We suggest that the critical probability $p_c \in [0,1]$ should be the maximizer of

$$\frac{\chi'(p)}{\chi(p)}.$$

Intuitively (and keeping in mind Russo's formula), this $p_c$ is the one at which adding a random edge to $G_p$ has the maximal impact on the size of the component containing $v$ (with relation to its size).

QUESTION 1. Do the properties mentioned above hold for this suggested $p_c$?

QUESTION 2. For "mean-field" graphs, such as the Hamming cube or the $d$-dimensional discrete torus for large $d$, does this definition coincide with the previous definition of Borgs, Chayes, van der Hofstad, Slade and Spencer [7], who require $\chi(p_c) = \lambda|V|^{1/3}$ (where $|V|$ is the number of vertices in the graph and $\lambda$ is some small constant)?

**Acknowledgments.** The first author would like to thank Microsoft Research, where this research was conducted, for their kind hospitality.

## REFERENCES

[1] ALDOUS, D. (1997). Brownian excursions, critical random graphs and the multiplicative coalescent. *Ann. Probab.* **25** 812–854. MR1434128
[2] ALDOUS, D. and FILL, J. (2007). *Reversible Markov Chains and Random Walks on Graphs.* Available at http://www.stat.berkeley.edu/~aldous/RWG/book.html.
[3] ALON, N. and SPENCER, J. H. (2000). *The Probabilistic Method*, 2nd ed. Wiley, New York. MR1885388
[4] BARLOW, M. and KUMAGAI, T. (2006). Random walk on the incipient infinite cluster on trees. Available at http://www.arxiv.org/abs/math.PR/0503118. MR2247823
[5] BENJAMINI, I., KOZMA, G. and WORMALD, N. (2006). The mixing time of the giant component of a random graph. Available at http://www.arxiv.org/abs/math.PR/0610459.
[6] BOLLOBÁS, B., JANSON, S. and RIORDAN, O. (2007). The phase transition in inhomogeneous random graphs. *Random Structures Algoritms* **31** 3–122. MR2337396
[7] BORGS, C., CHAYES, J. T., VAN DER HOFSTAD, R., SLADE, G. and SPENCER, J. (2005). Random subgraphs of finite graphs: I. The scaling window under the triangle condition. *Random Structures Algoritms* **27** 137–184. MR2155704

[8] CHANDRA, A. K., RAGHAVAN, P., RUZZO, W. L. and SMOLENSKY, R. (1989). The electrical resistance of a graph captures its commute and cover times. *Proc. Twenty-First ACM Symposium on Theory of Computing* 574–586. ACM Press, New York.

[9] CHUNG, F. and LU, L. (2001). The diameter of random sparse graphs. *Adv. in Appl. Math.* **26** 257–279. MR1826308

[10] ERDŐS, P. and RÉNYI, A. (1960). On the evolution of random graphs. *Magyar Tud. Akad. Mat. Kutató Int. Kőzl.* **5** 17–61. MR0125031

[11] FERNHOLZ, D. and RAMACHANDRAN, V. (2004). The diameter of sparse random graphs. Available at http://www.cs.utexas.edu/~fernholz/diam.ps. MR2362640

[12] FOUNTOULAKIS, N. and REED, B. A. (2006). The evolution of the mixing rate. Preprint.

[13] VAN DER HOFSTAD, R. and LUCZAK, M. (2006). Random subgraphs of the 2D Hamming graph: The supercritical phase. Preprint.

[14] KOLCHIN, V. F. (1986). *Random Mappings*. Optimization Software, New York. MR0865130

[15] KESTEN, H., NEY, P. and SPITZER, F. (1966). The Galton–Watson process with mean one and finite variance. *Theory Probab. Appl.* **11** 579–611. MR0207052

[16] KOLMOGOROV, A. N. (1938). On the solution of a problem in biology. *Izv. NII Matem. Mekh. Tomskogo Univ.* **2** 7–12.

[17] ŁUCZAK, T. (1998). Random trees and random graphs. *Random Structures Algoritms* **13** 485–500. MR1662797

[18] ŁUCZAK, T., PITTEL, B. and WIERMAN, J. C. (1994). The structure of a random graph at the point of the phase transition. *Trans. Amer. Math. Soc.* **341** 721–748. MR1138950

[19] LYONS, R. (1992). Random walks, capacity and percolation on trees. *Ann. Probab.* **20** 2043–2088. MR1188053

[20] LYONS, R., PEMANTLE, R. and PERES, Y. (1995). Conceptual proofs of $L \log L$ criteria for mean behavior of branching processes. *Ann. Probab.* **23** 1125–1138. MR1349164

[21] NACHMIAS, A. and PERES, Y. (2005). The critical random graph, with martingales. *Israel J. Math.* To appear. Available at http://www.arxiv.org/abs/math.PR/0512201.

[22] NACHMIAS, A. and PERES, Y. (2006). Critical percolation on random regular graphs. Preprint. Available at http://www.arxiv.org/abs/0707.2839.

[23] NASH-WILLIAMS, C. ST. J. A. (1959). Random walk and electric currents in networks. *Proc. Cambridge Philos. Soc.* **55** 181–194. MR0124932

[24] PERES, Y. (1999). Probability on trees: An introductory climb. *École d'Été de Probabilités de Saint-Flour XXVII. Lecture Notes in Math.* **1717** 193–280. Springer, Berlin. MR1746302

[25] TETALI, P. (1991). Random walks and the effective resistance of networks. *J. Theoret. Probab.* **4** 101–109. MR1088395

DEPARTMENT OF MATHEMATICS
UNIVERSITY OF CALIFORNIA, BERKELEY
BERKELEY, CALIFORNIA 94720
USA
E-MAIL: asafnach@math.berkeley.edu

MICROSOFT RESEARCH
ONE MICROSOFT WAY
REDMOND, WASHINGTON 98052-6399
USA
E-MAIL: peres@stat.berkeley.edu